\documentclass[11pt,a4paper]{amsart}
\usepackage[T1]{fontenc}
\usepackage{lmodern}
\usepackage{amsmath,amssymb,amsthm,mathtools}
\usepackage{microtype}
\usepackage[margin=28mm,top=27mm,bottom=28mm]{geometry}
\usepackage{xcolor}
\usepackage[colorlinks=true,linkcolor=blue!45!black,citecolor=blue!45!black,urlcolor=blue!45!black]{hyperref}
\hypersetup{pdftitle={Dimension functions and Sylvester rank functions},pdfsubject={BH\_v17}}
\newcommand{\doi}[1]{\href{https://doi.org/#1}{\nolinkurl{doi:#1}}}
\allowdisplaybreaks[1]
\numberwithin{equation}{section}
\newtheorem{theorem}{Theorem}[section]
\newtheorem{proposition}[theorem]{Proposition}
\newtheorem{lemma}[theorem]{Lemma}
\newtheorem{corollary}[theorem]{Corollary}
\newtheorem{maintheorem}{Theorem}

\theoremstyle{definition}
\newtheorem{definition}[theorem]{Definition}
\theoremstyle{remark}
\newtheorem{remark}[theorem]{Remark}
\newcommand{\DF}{\operatorname{DF}}
\newcommand{\LDF}{\operatorname{LDF}}
\newcommand{\Her}{\operatorname{Her}}
\newcommand{\eps}{\varepsilon}
\newcommand{\mpair}[2]{\left(\,#1\,\middle|\,#2\,\right)}
\newcommand{\mclass}[2]{\left[\left(\,#1\,\middle|\,#2\,\right)\right]}
\title[Dimension functions and Sylvester rank functions]{The space of Dimension functions on a $C^*$-algebra}

\author{Zhichao Liu}
\address{School of Mathematical Sciences,
Dalian University of Technology,
Dalian, {\rm 116024}, China }
\email{zhichaoliu@dlut.edu.cn}

\author{Tao Zhang}
\address{School of Mathematical Sciences,
Dalian University of Technology,
Dalian, {\rm 116024}, China }
\email{taozhang@mail.dlut.edu.cn}

\subjclass[2020]{Primary 46L05; Secondary 46L80, 46A55, 16D10}
\keywords{ dimension function, lattice, Sylvester rank function, Choquet simplex}

\begin{document}
\begin{abstract}
In this paper, we prove that the space of normalized dimension functions on a
unital $C^*$-algebra is either empty or a Choquet simplex. As a
byproduct, we show that $W(A)$ is hereditary in
$\operatorname{Cu}(A)$ for every $C^*$-algebra $A$.
The argument combines an exact factorization for Cuntz comparison
with Birkhoff's decomposition of lattice valuations. Moreover, we
show that the space
of normalized rank functions on a unital
ring is either empty or a Choquet simplex.
\end{abstract}
\maketitle
\section{Introduction}
Dimension functions on $C^*$-algebras were introduced by Cuntz
\cite{Cuntz}. They are nonnegative real-valued functions on positive
matrices, additive under direct sums and monotone for Cuntz comparison.
Cuntz's construction was motivated by Goodearl and Handelman's
work on rank functions for regular rings
\cite{GH76}.
Blackadar and Handelman developed the theory of dimension
functions and quasitraces in \cite{BH}.

For a unital $C^*$-algebra $A$, let $\DF(A)$ be the space of
dimension functions normalized by $d(1_A)=1$, with the topology
of pointwise convergence. Blackadar and Handelman proved that
the quasitracial state space is a Choquet simplex
\cite[Theorem~II.4.4]{BH}. They asked whether $\DF(A)$ is also
a Choquet simplex and whether the lower semicontinuous dimension
functions are dense in $\DF(A)$. The simplex question is also
recorded in \cite[Conjecture~12.2]{GP}.

Several authors have studied the simplex question for particular
classes of $C^*$-algebras. Perera proved the assertion for unital
algebras of real rank zero and stable rank one by establishing interpolation in the
ordered Grothendieck group of the Cuntz semigroup
\cite[Corollary 4.4]{Perera}. For simple unital finite
algebras which are either exact and $\mathcal Z$-stable or
approximately homogeneous with slow dimension growth, Brown,
Perera, and Toms proved both Blackadar--Handelman conjectures
\cite[Theorem~B]{BPT}. The interpolation argument also applies
to certain continuous fields, as Antoine, Bosa, Perera, and
Petzka showed \cite[Theorem~4.1]{ABPP}, and to all unital
algebras of stable rank one by the work of Antoine, Perera,
Robert, and Thiel \cite[Theorem~4.1]{APRT}. For an arbitrary unital
$C^*$-algebra, we obtain the simplex property from the lattice order
on the real span of the dimension functions.

\begin{maintheorem}\label{thm:dimension-simplex}
For every unital $C^*$-algebra $A$, the space $\DF(A)$ is either
empty or a Choquet simplex.
\end{maintheorem}

The factorization for Cuntz comparison used in the proof also yields
the following property of the Cuntz semigroup. Write $W(A)$ for the
Cuntz classes represented by positive elements in matrix algebras over $A$, and
$\operatorname{Cu}(A)$ for those represented in
$A\otimes\mathcal K$; see Definition~\ref{def:cuntz-semigroups}.
An ordered submonoid is hereditary if it contains every element
below one of its elements.

\begin{maintheorem}\label{thm:hereditary}
For every $C^*$-algebra $A$, $W(A)$ is
hereditary in $\operatorname{Cu}(A)$.
\end{maintheorem}

Blackadar, Robert, Tikuisis, Toms, and Winter proved hereditariness
for a possibly nonunital $C^*$-algebra $A$ provided that all
projections in every quotient of $A\otimes\mathcal K$ are finite
and that the radius of comparison is finite relative to some strictly
positive element of $A$ \cite[Theorem~4.4.1]{BRTTW}.
They asked whether hereditariness can fail for some $C^*$-algebra
\cite[Question~4.4.3]{BRTTW}.
Theorem~\ref{thm:hereditary} answers that question negatively.
The factorization gives the more precise statement that, if
$a\in(A\otimes\mathcal K)_+$ is Cuntz subequivalent to
$b\in M_q(A)_+$, then its class has a representative in
$M_{2q}(A)_+$.

The algebraic part of the proof applies to Sylvester matrix rank
functions over arbitrary unital rings. Malcolmson introduced
matrix and module rank functions in his study of homomorphisms
to skew fields \cite{Malcolmson}. We use the real-valued
formulation in \cite[Definition~2.7 and Theorem~2.9]{HL}; see
also \cite{GoodearlRegular} for rank functions on regular rings.
For a unital ring $R$, let $\mathcal R(R)$ be the cone of
Sylvester matrix rank functions, without normalization, and write
\[
 \mathbb P(R)=\{\rho\in\mathcal R(R):\rho(1_R)=1\}.
\]
All rank functions take values in $[0,\infty)$, and
$\mathbb P(R)$ has the topology of pointwise convergence.
The differences $\rho-\sigma$, with $\rho,\sigma\in\mathcal R(R)$,
form the real vector space $\mathcal R(R)-\mathcal R(R)$ under
pointwise addition and scalar multiplication. The partial order on $\mathcal{R}(R) - \mathcal{R}(R)$ is defined by
\[
 f\leq g\quad\Longleftrightarrow\quad g-f\in\mathcal R(R).
\]
Restricted to $\mathcal{R}(R)$, this order coincides with the algebraic order.
Hung and Li's criterion for the existence of a normalized rank function
compares the identity matrices $1_n\in M_n(R)$ in the Malcolmson preorder
$\preceq_{\mathrm M}$ of Definition~\ref{def:malcolmson}
\cite[Proposition~3.2]{HL}.

\begin{maintheorem}\label{thm:ring}
Let $R$ be an associative unital ring. Then the
ordered vector space $\mathcal R(R)-\mathcal R(R)$ is a vector
lattice. Consequently, $\mathbb P(R)$ is either empty or a Choquet simplex. Otherwise
\end{maintheorem}

To prove Theorem~\ref{thm:ring}, we use the modular lattice of
positive primitive formulas over $R$, written as matrix pairs;
see \cite{Prest,HLI}.
A rank function determines increasing valuations on these
lattices. Birkhoff's decomposition theorem
\cite[Theorem~17]{Birkhoff} supplies positive
parts for valuations of bounded variation. We prove that this
operation preserves additivity under direct sums of matrix pairs.
It therefore yields positive parts of differences of rank
functions, which establishes the vector lattice structure.


For a unital $C^*$-algebra, every dimension function extends
canonically to a Sylvester matrix rank function by
$d(X)=d(|X|)$, where $|X|=(X^*X)^{1/2}$
\cite[Propositions~2.6 and~5.2]{HL}. The converse fails
\cite[Examples~5.4 and~5.5]{HL}. We prove that every rank-function
summand of a dimension function is again a dimension function. Thus
dimension functions form a hereditary subcone of $\mathcal R(A)$, and
the canonical image of $\DF(A)$ is a closed face of
$\mathbb P(A)$. This proves Theorem~\ref{thm:dimension-simplex}.

This paper is organized as follows. In Section~\ref{sec:prelim}, we recall some basic facts on the space of dimension functions. In Section~\ref{sec:factor}, we prove the factorization result and Theorem~\ref{thm:hereditary}. In Section~\ref{sec:matrix-pairs}, we develop the rank-function argument and prove Theorem~\ref{thm:ring}. Finally, in Section~\ref{sec:refinement}, we use the techniques in Section~\ref{sec:matrix-pairs} 
to prove Theorem~\ref{thm:dimension-simplex}.

\section{Preliminaries}\label{sec:prelim}
Throughout Sections~\ref{sec:prelim} and~\ref{sec:factor}, we use $A$ to denote a unital $C^*$-algebra, unless otherwise specified.

For $p,q\geq1$, we write $M_{p,q}(A)$ for the $p\times q$ matrices
over $A$, and write $1_n$ and $0_n$ for the identity and zero
matrices of size $n$. 
Let
$M_\infty(A)=\bigcup_{n\geq1}M_n(A)$, with the usual upper-left
corner inclusions. Let $\mathcal K$ be the algebra
of compact operators on a separable infinite-dimensional Hilbert
space. For $a\in M_p(A)_+$ and $\eps>0$, $(a-\eps)_+$ denotes
the positive part of $a-\eps 1_p$. We write $\Her(a)$
for the hereditary $C^*$-subalgebra generated by $a$ and $a^{\oplus n}$
for the direct sum of $n$ copies of $a$. Thus
\[
 \Her(a)=\overline{aM_p(A)a}\quad(a\in M_p(A)_+),
 \qquad a^{\oplus n}=\operatorname{diag}(a,\ldots,a).
\]

\begin{definition}[{\cite[Definition~2.1]{GP}}]\label{def:cuntz}
Let $a\in M_p(A)_+$ and $b\in M_q(A)_+$. We say that $a$ is
\emph{Cuntz subequivalent} to $b$, and write $a\precsim b$, if
there are $x_j\in M_{q,p}(A)$ such that
\[
 \lim_{j\to\infty}\|x_j^*bx_j-a\|=0.
\]
We say that $a$ and $b$ are \emph{Cuntz equivalent}, and write
$a\sim b$, if $a\precsim b$ and $b\precsim a$.
The same definition applies to positive elements of $A\otimes\mathcal K$.
\end{definition}

The equivalences in the following proposition are the rectangular
matrix form of R{\o}rdam's comparison lemma \cite[Proposition~2.4]{Rordam}.

\begin{proposition}\label{prop:cuntz-basic}
Let $A$ be a $C^*$-algebra, $a\in M_p(A)_+$, and $b\in M_q(A)_+$.
The following are equivalent.
\begin{enumerate}
\item $a\precsim b$.
\item $(a-\eps)_+\precsim b$ for every $\eps>0$.
\item For every $\eps>0$ there are $\delta>0$ and $x\in M_{q,p}(A)$
such that $x^*x=(a-\eps)_+$ and $xx^*\in\Her((b-\delta)_+)$.
\end{enumerate}
Moreover, $X^*X\sim XX^*$ for $X\in M_{q,p}(A)$, and
$a^t\sim a$ for $t>0$.
\end{proposition}

\begin{proposition}\label{prop:hereditary-basic}
Let $A$ be a $C^*$-algebra and $b\in M_q(A)_+$.
\begin{enumerate}
\item If $c\in\Her(b)_+$, then $c\precsim b$ and $c\sim bcb$.
\item For $n\geq1$, $\Her(b^{\oplus n})=M_n(\Her(b))$.
If $x,y\in M_{q,p}(A)$ and $xx^*,yy^*\in\Her(b)$,
then $xy^*\in\Her(b)$.
\end{enumerate}
\end{proposition}
\begin{proof}
(1) For $\eps>0$, let $r_\eps=b(b+\eps1_q)^{-2}$ and
$e_\eps=br_\eps$, using the unitization when necessary.
Then $e_\eps$ is a positive contraction and
$\|(1_q-e_\eps)b\|\leq 2\eps$. Since $\Her(b)=\overline{bM_q(A)b}$,
$(e_\eps)$ is an approximate unit for $\Her(b)$.
For $c\in\Her(b)_+$, let $x_\eps=r_\eps^{1/2}c^{1/2}\in M_q(A)$.
Since $c^{1/2}\in\Her(b)$, we have
\[
 x_\eps^* b x_\eps=c^{1/2}e_\eps c^{1/2}\longrightarrow c,
 \qquad
 r_\eps(bcb)r_\eps=e_\eps c e_\eps\longrightarrow c.
\]
Thus $c\precsim b$ and $c\precsim bcb$. Since $bcb=b^*cb$, we also have $bcb\precsim c$. Hence $c\sim bcb$.

(2) Taking closures in
$b^{\oplus n}M_{nq}(A)b^{\oplus n}=M_n(bM_q(A)b)$ gives
$\Her(b^{\oplus n})=M_n(\Her(b))$. If $xx^*,yy^*\in\Her(b)$, then
\[
 \|x-e_\eps x\|^2
 =\|(1_q-e_\eps)xx^*(1_q-e_\eps)\|\longrightarrow0,
\]
and similarly $e_\eps y\to y$. It follows that
\[
 xy^*=\lim_{\eps\to0}e_\eps xy^*e_\eps
     =\lim_{\eps\to0}b(r_\eps xy^*r_\eps)b
     \in\Her(b). \qedhere
\]
\end{proof}

We use the positive-element formulation of dimension functions
in \cite[Definition~12.1]{GP}; compare
\cite[Definition~I.1.2]{BH}.

\begin{definition}\label{def:dimension}
A \emph{dimension function} on $A$ is a map
$d:M_\infty(A)_+\to[0,\infty)$ such that
\[
 d(0)=0,\qquad d(a\oplus b)=d(a)+d(b),\qquad
 a\precsim b\ \Longrightarrow\ d(a)\leq d(b).
\]
It is \emph{normalized} if $d(1_A)=1$.
Denote the set of normalized dimension functions on $A$ by
$\DF(A)$, equipped with the topology of pointwise convergence.
A dimension function $d$ is \emph{lower semicontinuous} if
\[
 d(a)=\sup_{\eps>0}d((a-\eps)_+)
 \qquad(a\in M_\infty(A)_+).
\]
Denote the set of normalized lower semicontinuous dimension
functions by $\LDF(A)$. This is the usual lower semicontinuity
condition; see \cite[p.~195]{BPT} and
\cite[Proposition~I.1.5]{BH}.
\end{definition}

Since $a\precsim1_p$ for $a\in M_p(A)_+$, every dimension
function satisfies
\[
 0\leq d(a)\leq p\,d(1_A).
\]
In particular, $d(1_A)=0$ implies $d=0$.
The defining conditions ensure that $\mathrm{DF}(A)$ forms a closed convex subset of the product space \(
\prod\limits_{p\geq1}\prod\limits_{a\in M_p(A)_+}[0,p]
\), subject to compatibility with corner inclusions. By Tychonoff's theorem, $\mathrm{DF}(A)$ is compact.

For arbitrary $C^*$-algebras, we use the Cuntz semigroups
introduced in \cite{Cuntz,CEI}; see also
\cite{GP}.

\begin{definition}\label{def:cuntz-semigroups}
Let $A$ be a $C^*$-algebra. Define
\[
 W(A)=M_\infty(A)_+/{\sim},\qquad
 \operatorname{Cu}(A)=(A\otimes\mathcal K)_+/{\sim}.
\]
We identify $W(A)$ with its canonical image in $\operatorname{Cu}(A)$. Let $[a]$ denote the Cuntz equivalence class of a positive element $a$. The natural partial order on $\operatorname{Cu}(A)$ is defined by $[a] \le [b]$ whenever $a \precsim b$, while addition is given by $[a] + [b] = [a \oplus b]$.
Following \cite{BRTTW},
we call this image \emph{hereditary} if
\[
 x\in\operatorname{Cu}(A),\quad y\in W(A),\quad x\leq y
 \quad\Longrightarrow\quad x\in W(A).
\]
\end{definition}

The following definition of a vector lattice is due to Kantorovich.
Freudenthal obtained it independently, drawing on the work of
F.~Riesz \cite[p.~238]{Birkhoff}.

\begin{definition}\label{def:ordered-space}
A \emph{real vector space} $E$ is a module over $\mathbb R$.
A partial order on $E$ makes it an \emph{ordered vector space} if $x\leq y$ implies
$x+z\leq y+z$ and $tx\leq ty$ for all $z\in E$ and $t\geq0$.
It is a \emph{vector lattice}, or \emph{Riesz space}, if every pair
$x,y\in E$ has a least upper bound $x\vee y$ and a greatest lower
bound $x\wedge y$.
\end{definition}

For a convex cone $C$ in a real vector space with $C\cap(-C)=\{0\}$,
its linear span is $C-C=\{x-y:x,y\in C\}$, with the vector
operations inherited from the ambient space. Its order is given by
$u\leq v$ if $v-u\in C$. For cones of nonnegative functions, this
order implies pointwise comparison and can be strictly stronger.

Let $K$ be a nonempty compact base of a convex cone $C$ in a locally
convex Hausdorff real vector space. By the criterion of Choquet and
Meyer \cite[Theorem~11]{CM}, $K$ is a Choquet simplex precisely when
$C-C$ is a vector lattice. 
For a cone with a compact base, the lattice condition is also equivalent
to interpolation and to refinement of finite sums in $C$.
The following proposition collects these characterizations \cite[Propositions~II.3.1--II.3.3
and~II.3.6]{Alfsen}, \cite[Proposition~2.1 and Theorem~2.2]{BR}.

\begin{proposition}\label{prop:simplex}
Let $C$ be a convex cone in a locally convex Hausdorff real vector space,
and let $\ell$ be a continuous linear functional which is strictly
positive on $C\setminus\{0\}$. Suppose
\[
 K=\{x\in C:\ell(x)=1\}
\]
is nonempty and compact. 
Then the following conditions are equivalent.
\begin{enumerate}
\item $K$ is a Choquet simplex.
\item $C-C$, ordered by $C$, is a vector lattice.
\item Whenever $a_1,a_2,b_1,b_2\in C$ satisfy $a_i\leq b_j$ for
all $i,j\in\{1,2\}$, there is a $c\in C$ such that
$a_i\leq c\leq b_j$ for all $i,j\in\{1,2\}$.
\item For all positive integers $r,s$ and all
$a_1,\ldots,a_r,b_1,\ldots,b_s\in C$ with
$\sum_{i=1}^r a_i=\sum_{j=1}^s b_j$, there are $c_{ij}\in C$ such that
\[
 a_i=\sum_{j=1}^s c_{ij}\quad(1\leq i\leq r),\qquad
 b_j=\sum_{i=1}^r c_{ij}\quad(1\leq j\leq s).
\]
\end{enumerate}
\end{proposition}

\section{A factorization for Cuntz comparison}\label{sec:factor}
We first assume that $A$ is unital and recall the following well-known
consequence of polynomial approximation.

\begin{proposition}\label{prop:intertwining}
Let $a\in M_p(A)_+$, $b\in M_q(A)_+$, and $x\in M_{q,p}(A)$.
If $xa=bx$, then $xf(a)=f(b)x$ for every continuous function $f$
on $[0,\max\{\|a\|,\|b\|\}]$. In particular, $xa^{1/2}=b^{1/2}x$.
\end{proposition}

\begin{proposition}\label{prop:factor}
Let $a\in M_p(A)_+$ and $b\in M_q(A)_+$ satisfy $a\precsim b$.
For every $n\geq1$, there is $T\in M_{(n+1)q,np}(A)$ such that
\[
 T^*T=a^{\oplus n},\qquad
 TT^*\in\Her(b^{\oplus(n+1)}).
\]
\end{proposition}

\begin{proof}
If $a=0$, take $T=0$. Suppose that $a\ne0$.
Choose $\eps_0>\|a\|$ and let $\eps_k=4^{-k}\eps_0$ for $k\geq0$.
Write $a_k=(a-\eps_k)_+$ and $B=\Her(b)$.
By Proposition~\ref{prop:cuntz-basic}(3), there are $x_k\in M_{q,p}(A)$ such that
\[
 x_k^*x_k=a_k,\qquad x_kx_k^*\in B.
\]
Let $x_k=u_ka_k^{1/2}$ be the polar decomposition in
$M_{q,p}(A^{**})$.

For $k\geq1$ and $j\in\{k,k+1\}$, define $g_{k,j}$ on
$[0,\|a\|]$ by
\[
 g_{k,j}(t)=
 \begin{cases}
  0,& t\leq\eps_{k-1},\\[2pt]
  \displaystyle\sqrt{\frac{t-\eps_{k-1}}{t-\eps_j}},&t>\eps_{k-1}.
 \end{cases}
\]
Since $\eps_j<\eps_{k-1}$, the function $g_{k,j}$ is continuous.
Thus $u_ja_{k-1}^{1/2}=x_jg_{k,j}(a)\in M_{q,p}(A)$, and
\begin{equation}\label{eq:cutoff-estimate}
 \|x_j-u_ja_{k-1}^{1/2}\|
 \leq\|a_j^{1/2}-a_{k-1}^{1/2}\|
 \leq\sqrt{\eps_{k-1}}
 \qquad(j=k,k+1).
\end{equation}

Choose $h_k=f_k(a)$, where $0\leq f_k\leq1$,
$f_k\in C_c((\eps_k,\infty))$, and $f_k=1$ on
$[\eps_{k-1},\|a\|]$. Define
\[
 v_k=u_kh_ku_{k+1}^*
     =(u_kh_k^{1/2})(u_{k+1}h_k^{1/2})^*.
\]
By functional calculus, both factors on the right are matrices
over $A$, and their products with their adjoints belong to $B$.
Thus $v_k\in B$ by Proposition~\ref{prop:hereditary-basic}(2), and
$\|v_k\|\leq1$.

Since $h_ka_{k-1}=a_{k-1}$, the support projections in the polar
decompositions give
\[
 v_ku_{k+1}a_{k-1}^{1/2}=u_ka_{k-1}^{1/2},\qquad
 v_k^*u_ka_{k-1}^{1/2}=u_{k+1}a_{k-1}^{1/2}.
\]
It follows that $1_q-v_k^*v_k$ and its positive square root both
annihilate $u_{k+1}a_{k-1}^{1/2}$. Applying
Proposition~\ref{prop:intertwining} to $v_k(1_q-v_k^*v_k)=(1_q-v_kv_k^*)v_k$ gives
\[
 v_k(1_q-v_k^*v_k)^{1/2}=(1_q-v_kv_k^*)^{1/2}v_k.
\]
Using this equality and its adjoint, we see that the matrix
\[
 U_k=\begin{pmatrix}
 v_k&(1_q-v_kv_k^*)^{1/2}\\
 (1_q-v_k^*v_k)^{1/2}&-v_k^*
 \end{pmatrix}\in M_2(B+\mathbb C1_q)
\]
is unitary and satisfies
\begin{equation}\label{eq:unitary}
 U_k\binom{u_{k+1}a_{k-1}^{1/2}}{0}
 =\binom{u_ka_{k-1}^{1/2}}{0}.
\end{equation}

Fix $n\geq1$. For $1\leq j\leq n$, let $U_{k,j}$ act as $U_k$
on the $j$th and $(n+1)$st block rows and as the identity on the
other block rows. With block sizes $(j-1)q,q,(n-j)q,q$, this is
\[
 U_{k,j}=\begin{pmatrix}
 1_{(j-1)q}&0&0&0\\
 0&v_k&0&(1_q-v_kv_k^*)^{1/2}\\
 0&0&1_{(n-j)q}&0\\
 0&(1_q-v_k^*v_k)^{1/2}&0&-v_k^*
 \end{pmatrix}.
\]
Blocks of size zero are omitted. Define
\[
 W_k=U_{k,n}U_{k,n-1}\cdots U_{k,1}
 \in M_{n+1}(B+\mathbb C1_q)\subseteq M_{(n+1)q}(A).
\]
By \eqref{eq:unitary}, each successive factor changes one diagonal
entry from $u_{k+1}a_{k-1}^{1/2}$ to $u_ka_{k-1}^{1/2}$ and leaves
the last block row zero. Induction on the block rows gives
\begin{equation}\label{eq:W-action}
 W_k\begin{pmatrix}
 u_{k+1}a_{k-1}^{1/2}&&0\\
 &\ddots&\\
 0&&u_{k+1}a_{k-1}^{1/2}\\
 0&\cdots&0
 \end{pmatrix}
 =\begin{pmatrix}
 u_ka_{k-1}^{1/2}&&0\\
 &\ddots&\\
 0&&u_ka_{k-1}^{1/2}\\
 0&\cdots&0
 \end{pmatrix}.
\end{equation}

In $M_{n+1}(B+\mathbb C1_q)$ define $V_1=1_{(n+1)q}$ and
$V_{k+1}=V_kW_k$. Let
\[
 T_k=V_k\begin{pmatrix}
 x_k&&0\\ &\ddots&\\0&&x_k\\0&\cdots&0
 \end{pmatrix}.
\]
Since $V_k$ and $W_k$ are unitary, equations
\eqref{eq:cutoff-estimate} and \eqref{eq:W-action} give
\[
 \begin{aligned}
 \|T_{k+1}-T_k\|
 &=\left\|
 W_k\begin{pmatrix}x_{k+1}&&0\\&\ddots&\\0&&x_{k+1}\\0&\cdots&0\end{pmatrix}
 -\begin{pmatrix}x_k&&0\\&\ddots&\\0&&x_k\\0&\cdots&0\end{pmatrix}
 \right\|\\
 &\leq\|x_{k+1}-u_{k+1}a_{k-1}^{1/2}\|
       +\|x_k-u_ka_{k-1}^{1/2}\|\\
 &\leq2\sqrt{\eps_{k-1}}.
 \end{aligned}
\]
Since $\sum_{k=1}^{\infty}\sqrt{\eps_{k-1}}<\infty$,
the sequence $(T_k)$ converges in norm to a matrix
$T\in M_{(n+1)q,np}(A)$.
Since $T_k^*T_k=a_k^{\oplus n}$ and $a_k\to a$, we obtain
$T^*T=a^{\oplus n}$.
Moreover, $M_{n+1}(B)$ is an ideal in
$M_{n+1}(B+\mathbb C1_q)$, so
\[
 T_kT_k^*
 =V_k\operatorname{diag}(x_kx_k^*,\ldots,x_kx_k^*,0_q)V_k^*
 \in M_{n+1}(B).
\]
Taking the norm limit gives
$TT^*\in M_{n+1}(B)=\Her(b^{\oplus(n+1)})$.
\end{proof}

The case $n=1$ gives the following more precise form of
Theorem~\ref{thm:hereditary}. Here $A$ need not be unital.

\begin{corollary}\label{cor:hereditary}
Let $A$ be a $C^*$-algebra, not necessarily unital.
If
$a\in(A\otimes\mathcal K)_+$, $b\in M_q(A)_+$, and $a\precsim b$,
then $a\sim c$ for some $c\in M_{2q}(A)_+$.
If $A$ is unital, then
\[
 W(A)=\{x\in\operatorname{Cu}(A):x\leq m[1_A]
              \text{ for some }m\geq1\}.
\]
\end{corollary}

\begin{proof}
Let $D=A\otimes\mathcal K$ and regard $M_q(A)$ as the usual
finite corner in $D$. Apply Proposition~\ref{prop:factor} in
the unitization $\widetilde D$ to $a,b\in\widetilde D_+$, with
$n=1$. This gives a
column $T=\binom{t_1}{t_2}$ over $\widetilde D$ such that
\[
 T^*T=a,\qquad
 TT^*=\begin{pmatrix}t_1t_1^*&t_1t_2^*\\
                       t_2t_1^*&t_2t_2^*\end{pmatrix}
 \in\Her_{M_2(\widetilde D)}(b\oplus b).
\]
For the quotient map $\pi:\widetilde D\to\widetilde D/D$,
the equality $T^*T=a\in D$ gives $\pi(T)^*\pi(T)=0$.
Thus $T\in M_{2,1}(D)$. Since $b$ lies in the corner $M_q(A)$,
\[
 \Her_{M_2(\widetilde D)}(b\oplus b)
 =M_2(\Her_D(b))\subseteq M_{2q}(A).
\]
Consequently, $c=TT^*$ belongs to $M_{2q}(A)_+$ and
$a=T^*T\sim TT^*=c$ in the stabilization of $D$.
This proves the first assertion and Theorem~\ref{thm:hereditary}.
If $A$ is unital, every
$b\in M_q(A)_+$ satisfies $[b]\leq q[1_A]$.
The converse inclusion in the displayed formula now follows
from hereditariness.
\end{proof}

\section{Sylvester rank functions on unital rings}\label{sec:matrix-pairs}
In this section we prove Theorem~\ref{thm:ring}. Throughout,
$R$ is an associative unital ring, and modules are unital left
$R$-modules. The notation $M_{p,q}(R)$, $1_n$, and $0_n$ has
the same meaning as in Section~\ref{sec:prelim}. Matrices with zero rows or zero columns are allowed.

\begin{definition}[{\cite[Definition~2.7]{HL}}]\label{def:rank}
A \emph{Sylvester matrix rank function} on $R$ is a function $\rho$
from all rectangular matrices over $R$ to $[0,\infty)$ such that
\[
 \begin{gathered}
 \rho(0)=0,\qquad
 \rho(XY)\leq\min\{\rho(X),\rho(Y)\},\\
 \rho\!\left(\begin{pmatrix}X&0\\0&Y\end{pmatrix}\right)=\rho(X)+\rho(Y),\qquad
 \rho\!\left(\begin{pmatrix}X&Z\\0&Y\end{pmatrix}\right)\geq\rho(X)+\rho(Y).
 \end{gathered}
\]
It is \emph{normalized} if $\rho(1_R)=1$. Denote the cone of these
functions by $\mathcal R(R)$ and its normalized subset by $\mathbb P(R)$.
Matrices with zero rows or zero columns have rank zero.
\end{definition}

As with dimension functions, normalization is omitted from the
definition.

\begin{definition}[{\cite[Definition~3.1]{HL}}]\label{def:malcolmson}
Let $X$ and $Y$ be rectangular matrices over $R$.
We write $X\precsim_{\mathrm M}Y$ if either $X=CYD$ for compatible
matrices $C,D$, or
\[
 Y=\begin{pmatrix}C&E\\0&D\end{pmatrix},\qquad
 X=\begin{pmatrix}C&0\\0&D\end{pmatrix}
\]
for compatible matrices $C,D,E$.
We write $X\preceq_{\mathrm M}Y$ if there is a finite chain
\[
 X=X_1\precsim_{\mathrm M}X_2\precsim_{\mathrm M}\cdots
 \precsim_{\mathrm M}X_m=Y.
\]
The relation $\preceq_{\mathrm M}$ is called \emph{Malcolmson subequivalence}.
We write $X\sim_{\mathrm M}Y$ if $X\preceq_{\mathrm M}Y$ and
$Y\preceq_{\mathrm M}X$.
\end{definition}

Thus $\precsim_{\mathrm M}$ denotes the elementary relation and
$\preceq_{\mathrm M}$ its transitive closure. 
Then we record the elementary consequences of the rank axioms used below.

\begin{proposition}\label{prop:rank-properties}
Let $\rho\in\mathcal R(R)$. Then we have  the following.
\begin{enumerate}
\item For $X\in M_{k,n}(R)$,
\[
 0\leq\rho(X)\leq\min\{k,n\}\rho(1_R).
\]
If $\rho(1_R)=0$, then $\rho=0$.
\item For $X\in M_{k,n}(R)$ and invertible matrices $U\in M_k(R)$,
$V\in M_n(R)$, one has $\rho(UXV)=\rho(X)$.
\item Adjoining zero rows or zero columns to a matrix preserves its
rank. For $X\in M_{k,n}(R)$ one also has
\[
 \rho\!\left(\begin{pmatrix}1_k&X\end{pmatrix}\right)=\rho(1_k).
\]
\item If $X$ and $Y$ have the same number of rows, then
\[
 \max\{\rho(X),\rho(Y)\}
 \leq\rho\!\left(\begin{pmatrix}X&Y\end{pmatrix}\right)
 \leq\rho(X)+\rho(Y).
\]
\item If $X$ and $Y$ have the same number of columns, then
\[
 \max\{\rho(X),\rho(Y)\}
 \leq\rho\!\left(\binom{X}{Y}\right)
 \leq\rho(X)+\rho(Y).
\]
\item If $X$ and $Y$ have the same size, then
$\rho(X+Y)\leq\rho(X)+\rho(Y)$.
\item If $X\in M_{k,n}(R)$ and $Y\in M_{n,m}(R)$, then
\[
 \rho(XY)\geq\rho(X)+\rho(Y)-\rho(1_n).
\]
\item If $X\preceq_{\mathrm M}Y$, then $\rho(X)\leq\rho(Y)$.
\end{enumerate}
\end{proposition}

\begin{proof}
Since $\rho(1_k)=k\rho(1_R)$ and $X=1_kX=X1_n$, the product
inequality gives (1). For (2), the same inequality gives
\[
 \rho(UXV)\leq \rho(X)
   =\rho\bigl(U^{-1}(UXV)V^{-1}\bigr)\leq \rho(UXV).
\]
For (3), the two factorizations
\[
 \begin{pmatrix}X&0\\0&0\end{pmatrix}
 =\binom{1_k}{0}X\begin{pmatrix}1_n&0\end{pmatrix},
 \]\[
 X=\begin{pmatrix}1_k&0\end{pmatrix}
       \begin{pmatrix}X&0\\0&0\end{pmatrix}\binom{1_n}{0}
\]
give both rank inequalities. Either added block may have size zero.
Moreover,
\[
 \begin{pmatrix}1_k&X\end{pmatrix}
 =\begin{pmatrix}1_k&0\end{pmatrix}
   \begin{pmatrix}1_k&X\\0&1_n\end{pmatrix}.
\]
The square factor is invertible, so (2) and the zero-column equality
give the stated rank.

For (4), multiplying $\begin{pmatrix}X&Y\end{pmatrix}$ by the
coordinate inclusions gives $X$ and $Y$, hence the lower bound.
If the two matrices have $k$ rows, the upper bound follows from
\[
 \begin{pmatrix}X&Y\end{pmatrix}
 =\begin{pmatrix}1_k&1_k\end{pmatrix}(X\oplus Y).
\]
For (5), the coordinate projections recover $X$ and $Y$ from
$\binom{X}{Y}$. If they have $n$ columns, then
\[
 \binom{X}{Y}=(X\oplus Y)\binom{1_n}{1_n}.
\]
The product inequality and additivity on block sums give the two bounds.
For (6), let $X$ and $Y$ have $n$ columns. By the product inequality
and (4),
\[
 \rho(X+Y)
 =\rho\!\left(\begin{pmatrix}X&Y\end{pmatrix}\binom{1_n}{1_n}\right)
 \leq\rho\!\left(\begin{pmatrix}X&Y\end{pmatrix}\right)
 \leq\rho(X)+\rho(Y).
\]

For (7), elementary block operations give
\[
 \begin{pmatrix}1_n&0\\-X&1_k\end{pmatrix}
 \begin{pmatrix}1_n&Y\\X&0\end{pmatrix}
 \begin{pmatrix}1_n&-Y\\0&1_m\end{pmatrix}
 =\begin{pmatrix}1_n&0\\0&-XY\end{pmatrix}.
\]
A column permutation changes the middle matrix to
$\begin{pmatrix}Y&1_n\\0&X\end{pmatrix}$.
Its rank is at least $\rho(Y)+\rho(X)$, proving the inequality.
Each elementary comparison in Definition~\ref{def:malcolmson}
preserves the rank inequality, so transitivity proves (8).
\end{proof}

We give $\mathbb P(R)$ the topology of pointwise convergence.
The rank axioms are closed conditions in this topology.
Proposition~\ref{prop:rank-properties}(1) therefore shows that
$\mathbb P(R)$ is a compact convex set. Every nonzero element of
$\mathcal R(R)$ has a unique expression $t\rho$ with $t>0$ and
$\rho\in\mathbb P(R)$.

On $\mathcal R(R)-\mathcal R(R)$, we use the order
\[
 f\leq g\quad\Longleftrightarrow\quad g-f\in\mathcal R(R).
\]
It restricts to the algebraic order on $\mathcal R(R)$.
We will construct the positive part of each difference of rank
functions by using valuations on the lattices of matrix pairs.

We use the standard lattice of positive primitive formulas over
$R$.

\begin{definition}[{\cite{Prest,HLI}}]\label{def:matrix-pairs}
Fix $n\geq1$. A \emph{positive primitive formula}, or
\emph{pp-formula}, in $n$ free variables over $R$ has the form
\[
 \exists\,y\ (By=Cx),\qquad
 B\in M_{k,l}(R),\quad C\in M_{k,n}(R),
\]
where $x$ and $y$ are columns of lengths $n$ and $l$.
We represent the formula by the \emph{matrix pair} $\mpair{B}{C}$.
Its solution set in a left $R$-module $M$ is
\[
 \mpair{B}{C}(M)
 =\{x\in M^n:By=Cx\text{ for some }y\in M^l\}.
\]
This is an additive subgroup of $M^n$, called a
\emph{pp-definable subgroup} or \emph{finite matrix subgroup}.
It need not be an $R$-submodule when $R$ is noncommutative.
The case $l=0$ is allowed and gives the equation $Cx=0$.

For two pairs with $n$ free variables, write
\[
 \mpair{B}{C}\leq_n\mpair{D}{E}
\]
if $\mpair{B}{C}(M)\subseteq\mpair{D}{E}(M)$ for every left
$R$-module $M$.
Write $\mpair{B}{C}\sim_n\mpair{D}{E}$ when both comparisons
hold. The class of $\mpair{B}{C}$ is denoted by
$\mclass{B}{C}$, and the ordered set of these classes by
$\mathcal L_n(R)$. If $P=\mclass{B}{C}$, we write
$P(M)=\mpair{B}{C}(M)$.
\end{definition}

For example, let $X\in M_{k,n}(R)$ and let $M$ be a left $R$-module.
The class $\mclass{0}{X}\in\mathcal L_n(R)$ is given by the equation
$Xx=0$, so its solution set is the kernel of the additive map
$M^n\to M^k$, $x\mapsto Xx$. For $\mclass{X}{1_k}\in\mathcal L_k(R)$
the free variable is $z\in M^k$. The condition on $z$ is $z=Xy$ for
some $y\in M^n$, so the solution set is the image of the same map.

The comparison criterion below is \cite[Lemma~1.1.13]{Prest},
written for left modules. The lattice structure and its operations
are described in \cite{HLI}.

\begin{proposition}\label{prop:pair-modular}
Let $B\in M_{k,l}(R)$, $C\in M_{k,n}(R)$,
$D\in M_{s,t}(R)$, and $E\in M_{s,n}(R)$.
\begin{enumerate}
\item The relation $\mpair{B}{C}\leq_n\mpair{D}{E}$ holds
if and only if there exist
\[
 U\in M_{s,k}(R),\qquad V\in M_{t,l}(R),\qquad W\in M_{t,n}(R)
\]
such that
\begin{equation*}
 UB=DV,\qquad UC=E+DW.
\end{equation*}
\item The ordered set $\mathcal L_n(R)$ is a bounded modular
lattice. For $P=\mclass{B}{C}$ and $Q=\mclass{D}{E}$,
\begin{equation*}\label{eq:pair-lattice}
 \begin{aligned}
 P\wedge Q&=\mclass{\begin{pmatrix}B&0\\0&D\end{pmatrix}}
                      {\binom{C}{E}},\\
 P\vee Q&=\mclass{\begin{pmatrix}C&B&0\\-E&0&D\end{pmatrix}}
                    {\binom{0}{E}}.
 \end{aligned}
\end{equation*}
For each left $R$-module $M$, we have
\[
 (P\wedge Q)(M)=P(M)\cap Q(M),\qquad
 (P\vee Q)(M)=P(M)+Q(M),
\]
where $P(M)+Q(M)$ is the additive subgroup generated by $P(M)\cup Q(M)$.
The least and greatest elements are $\mclass{0}{1_n}$ and
$\mclass{0}{0_n}$, respectively.
\end{enumerate}
\end{proposition}

For pairs with $n$ and $m$ coordinates, respectively, set
\[
 \mpair{B}{C}\oplus\mpair{D}{E}
 =\mpair{B\oplus D}{C\oplus E}.
\]
This operation is well defined on equivalence classes because
\begin{equation*}
 (P\oplus Q)(M)=P(M)\times Q(M)\subseteq M^{n+m}.
\end{equation*}
It preserves order in each variable. In particular,
\[
 \begin{aligned}
 (P_1\oplus Q_1)\wedge(P_2\oplus Q_2)
   &=(P_1\wedge P_2)\oplus(Q_1\wedge Q_2),\\
 (P_1\oplus Q_1)\vee(P_2\oplus Q_2)
   &=(P_1\vee P_2)\oplus(Q_1\vee Q_2).
 \end{aligned}
\]
Here $P_i\in\mathcal L_n(R)$ and $Q_i\in\mathcal L_m(R)$.
We use $\wedge$ and $\vee$ for the lattice operations on classes.

We next recall Birkhoff's theorem for valuations on a modular
lattice.

\begin{definition}[{\cite{Birkhoff}}]\label{def:valuations}
Let $\mathcal L$ be a bounded lattice with least element $0_{\mathcal L}$ and greatest
element $1_{\mathcal L}$. A \emph{valuation} is a function $v:\mathcal L\to\mathbb R$
such that
\[
 v(x)+v(y)=v(x\wedge y)+v(x\vee y).
\]
It is \emph{increasing} if $x\leq y$ implies $v(x)\leq v(y)$.
It has \emph{bounded variation} if
\[
 \sup_{0_{\mathcal L}=x_0\leq\cdots\leq x_r=1_{\mathcal L}}
 \sum_{j=1}^r|v(x_j)-v(x_{j-1})|<\infty.
\]
On valuations vanishing at $0_{\mathcal L}$, define the \emph{valuation order}
by
\[
 v\leq w\quad\Longleftrightarrow\quad w-v\text{ is increasing}.
\]
Equivalently, $v(y)-v(x)\leq w(y)-w(x)$ whenever $x\leq y$.
This order is stronger than pointwise order.
\end{definition}

\begin{theorem}[Birkhoff {\cite[Theorem~17]{Birkhoff}}]\label{lem:birkhoff}
Let $\mathcal L$ be a bounded modular lattice. Its real-valued valuations
of bounded variation which vanish at $0_{\mathcal L}$ form a vector lattice
in the valuation order. For such a valuation $u$, the positive part
$u^+=u\vee0$ is
\begin{equation*}
 u^+(x)=\sup_{0_{\mathcal L}=x_0\leq\cdots\leq x_r=x}
       \sum_{j=1}^r\max\{u(x_j)-u(x_{j-1}),0\}.
\end{equation*}
The supremum is over finite chains; the empty chain has value zero.
In particular, $u^+$ is an increasing valuation and is the least
upper bound of $u$ and zero in the valuation order.
\end{theorem}

We have taken Birkhoff's fixed base point to be $0_{\mathcal L}$;
see
\cite[pp.~83--84, formulas~(11)--(14)]{Birkhoff}.
His term ``isotone'' means increasing.
Compatibility with direct sums will be proved in
Lemma~\ref{lem:positive-block}.

\begin{definition}\label{def:valuation-family}
Let $v=(v_n)_{n\geq1}$, where $v_n:\mathcal L_n(R)\to\mathbb R$
is a valuation for each $n$. We write $v(P)=v_n(P)$ when
$P\in\mathcal L_n(R)$. The family is \emph{increasing} if each $v_n$
is increasing. It \emph{vanishes at the origin} if
$v_n(\mclass{0}{1_n})=0$ for every $n$, and is \emph{additive} if
\[
 v_{n+m}(P\oplus Q)=v_n(P)+v_m(Q)
 \qquad(P\in\mathcal L_n(R),\ Q\in\mathcal L_m(R),\ n,m\geq1).
\]
Families vanishing at the origin are ordered componentwise by the
valuation order.
\end{definition}

\begin{lemma}\label{lem:rank-valuation}
Let $\rho\in\mathcal R(R)$. For $\mclass{B}{C}\in\mathcal L_n(R)$, define
\[
 v_\rho(\mclass{B}{C})=n\,\rho(1_R)-\rho\!\left(\begin{pmatrix}C&B\end{pmatrix}\right)+\rho(B).
\]
Then $v_\rho$ is a well-defined increasing additive family of
valuations vanishing at the origin. Moreover,
\[
 0\leq v_\rho(P)\leq n\,\rho(1_R),\qquad
 v_\rho(\mclass{X}{1_k})=\rho(X)\quad(X\in M_{k,n}(R)).
\]
\end{lemma}

\begin{proof}
Let $P=\mclass{B}{C}$ and $Q=\mclass{D}{E}$ have the matrix
sizes in Proposition~\ref{prop:pair-modular}, and set
\[
 L=\begin{pmatrix}C&B&0\\-E&0&D\end{pmatrix}.
\]
Suppose that $P\leq_n Q$, and choose $U,V,W$ as in Proposition \ref{prop:pair-modular}(1). Since
\[
 \begin{pmatrix}1_k&0\\U&1_s\end{pmatrix}
 L\begin{pmatrix}1_n&0&0\\0&1_l&0\\-W&-V&1_t\end{pmatrix}
 =\begin{pmatrix}C&B&0\\0&0&D\end{pmatrix},
\]
Proposition~\ref{prop:rank-properties}(2) and block-sum additivity
give $$\rho(L)=\rho\!\left(\begin{pmatrix}C&B\end{pmatrix}\right)+\rho(D).$$
On the other hand, we also have
\[
 \rho(L)=\rho\!\left(\begin{pmatrix}B&C&0\\0&-E&D\end{pmatrix}\right)
 \geq \rho(B)+\rho\!\left(\begin{pmatrix}-E&D\end{pmatrix}\right)
 =\rho(B)+\rho\!\left(\begin{pmatrix}E&D\end{pmatrix}\right).
\]
It follows that
\[
 \begin{aligned}
 v_\rho(Q)-v_\rho(P)
 &=\rho\!\left(\begin{pmatrix}C&B\end{pmatrix}\right)+\rho(D)
   -\rho\!\left(\begin{pmatrix}E&D\end{pmatrix}\right)-\rho(B)\\
 &=\rho(L)-\rho(B)-\rho\!\left(\begin{pmatrix}E&D\end{pmatrix}\right)\geq0.
 \end{aligned}
\]
Thus $v_\rho$ is increasing. Applying the inequality in both
directions shows that it is well defined on equivalence classes.

For arbitrary $P,Q$, the full coefficient matrix of $P\wedge Q$
satisfies
\[
 \begin{pmatrix}C&B&0\\E&0&D\end{pmatrix}
 =\operatorname{diag}(1_k,-1_s)\,L\,
    \operatorname{diag}(1_n,1_l,-1_t).
\]
Thus its rank is $\rho(L)$ by Proposition~\ref{prop:rank-properties}(2),
and
\[
 v_\rho(P\wedge Q)=n\,\rho(1_R)-\rho(L)+\rho(B)+\rho(D).
\]
Also,
\[
 P\vee Q=\mclass{\begin{pmatrix}C&B&0\\-E&0&D\end{pmatrix}}
                  {\begin{pmatrix}0\\E\end{pmatrix}}.
\]
Right multiplication of its full coefficient matrix by an invertible
matrix gives
\[
 \begin{pmatrix}0&C&B&0\\E&-E&0&D\end{pmatrix}
 \begin{pmatrix}1_n&1_n&0&0\\1_n&0&0&0\\0&0&1_l&0\\0&0&0&1_t\end{pmatrix}
 =\begin{pmatrix}C&0&B&0\\0&E&0&D\end{pmatrix}.
\]
After a column permutation, the matrix on the right is
$\begin{pmatrix}C&B\end{pmatrix}\oplus\begin{pmatrix}E&D\end{pmatrix}$.
Hence
\[
 v_\rho(P\vee Q)=n\,\rho(1_R)-\rho\!\left(\begin{pmatrix}C&B\end{pmatrix}\right)
 -\rho\!\left(\begin{pmatrix}E&D\end{pmatrix}\right)+\rho(L).
\]
Thus
\[
 \begin{aligned}
 v_\rho(P\wedge Q)+v_\rho(P\vee Q)
 &=2n\,\rho(1_R)-\rho\!\left(\begin{pmatrix}C&B\end{pmatrix}\right)
 -\rho\!\left(\begin{pmatrix}E&D\end{pmatrix}\right)+\rho(B)+\rho(D)\\
 &=v_\rho(P)+v_\rho(Q).
 \end{aligned}
\]

Let $P=\mclass{B}{C}\in\mathcal L_n(R)$ and $Q=\mclass{D}{E}\in\mathcal L_m(R)$.
A column permutation carries $(C\oplus E\;\;B\oplus D)$ to
$(C\;\;B)\oplus(E\;\;D)$, so block-sum additivity of $\rho$ gives
$$v_\rho(P\oplus Q)=v_\rho(P)+v_\rho(Q). $$

Finally, it is straightforward to verify that
$$v_\rho(\mclass{0}{0_n})=n\,\rho(1_R),\qquad 0\leq v_\rho(P)\leq n\,\rho(1_R).$$
For $X\in M_{k,n}(R)$, Proposition~\ref{prop:rank-properties}(3) gives
$\rho\bigl((1_k\;\;X)\bigr)=\rho(1_k)=k\rho(1_R)$, whence
$$v_\rho(\mclass{X}{1_k})=k\rho(1_R)-\rho(1_k)+\rho(X)=\rho(X).$$
\end{proof}

\begin{proposition}\label{prop:image-kernel}
Let $v$ be an additive family of valuations vanishing at the origin.
For $X\in M_{k,n}(R)$, with $k,n\geq1$, one has
\[
 v(\mclass{X}{1_k})+v(\mclass{0}{X})
 =v(\mclass{0}{0_n})=n\,v(\mclass{0}{0_1}).
\]
\end{proposition}

\begin{proof}
Consider the class
\[
 G_X=\mclass{0}{\begin{pmatrix}-X&1_k\end{pmatrix}}
 \quad\text{in }\mathcal L_{n+k}(R).
\]
For every left $R$-module $M$, its solution set is the graph
$\{(x,Xx):x\in M^n\}$. This graph meets $\{0\}\times M^k$ in
zero, and each $(x,z)\in M^{n+k}$ is the sum of $(x,Xx)$ and
$(0,z-Xx)$. The valuation identity for $G_X$ and
$\mclass{0}{1_n}\oplus\mclass{0}{0_k}$ therefore gives
\[
 v(G_X)+v(\mclass{0}{0_k})=v(\mclass{0}{0_{n+k}})
 =v(\mclass{0}{0_n})+v(\mclass{0}{0_k}).
\]
Hence $v(G_X)=v(\mclass{0}{0_n})$.

The intersection of the graph with $M^n\times\{0\}$ consists of
the pairs $(x,0)$ with $Xx=0$. The sum of the graph and $M^n\times\{0\}$ consists of
all pairs $(x,z)$ with $x\in M^n$ and $z=Xy$ for some $y\in M^n$,
since $(x,Xy)=(y,Xy)+(x-y,0)$. Consequently,
\[
 \begin{aligned}
 G_X\wedge\bigl(\mclass{0}{0_n}\oplus\mclass{0}{1_k}\bigr)
   &=\mclass{0}{X}\oplus\mclass{0}{1_k},\\
 G_X\vee\bigl(\mclass{0}{0_n}\oplus\mclass{0}{1_k}\bigr)
   &=\mclass{0}{0_n}\oplus\mclass{X}{1_k}.
 \end{aligned}
\]
Applying the valuation identity and cancelling $v(\mclass{0}{0_n})$
gives
\[
 v(\mclass{0}{X})+v(\mclass{X}{1_k})=v(G_X)
 =v(\mclass{0}{0_n})=n\,v(\mclass{0}{0_1}),
\]
where the last equality follows from additivity.
\end{proof}

\begin{corollary}\label{lem:valuation-matrix}
Let $v$ be an increasing additive family of valuations vanishing at
the origin. For $X\in M_{k,n}(R)$ with $k,n\geq1$, define
\[
 f(X)=v(\mclass{X}{1_k}).
\]
Then $f$ extends uniquely to an element of $\mathcal R(R)$, and
$f(1_R)=v(\mclass{0}{0_1})$.
\end{corollary}

\begin{proof}
Let $c=v(\mclass{0}{0_1})$. Since
$\mclass{0}{1_k}\leq_k\mclass{X}{1_k}\leq_k\mclass{0}{0_k}$,
we have $0\leq f(X)\leq kc$. The equality $f(0)=0$ and additivity
on block sums follow from Definition~\ref{def:valuation-family}.

If $X\in M_{k,n}(R)$ and $Y\in M_{n,m}(R)$, image inclusion gives
$\mclass{XY}{1_k}\leq_k\mclass{X}{1_k}$, so $f(XY)\leq f(X)$.
Kernel inclusion gives $\mclass{0}{Y}\leq_m\mclass{0}{XY}$.
Proposition~\ref{prop:image-kernel} now gives
\[
 f(XY)=mc-v(\mclass{0}{XY})\leq mc-v(\mclass{0}{Y})=f(Y).
\]
Thus the product inequality holds.

For the upper triangular inequality, let
\[
 T=\begin{pmatrix}X&Z\\0&Y\end{pmatrix},
 \quad X\in M_{k,n}(R),\quad Y\in M_{p,q}(R),\quad Z\in M_{k,q}(R),
\]
and define $P,Q\in\mathcal L_{k+p}(R)$ by
\[
 P=\mclass{T}{1_{k+p}},\qquad
 Q=\mclass{0}{0_k}\oplus\mclass{0}{1_p}.
\]
Then $f(T)=v(P)$. For every left module $M$, $P(M)$ is the image of $T$ and
$Q(M)=M^k\times\{0\}$. Consequently,
\[
 \mclass{X}{1_k}\oplus\mclass{0}{1_p}\leq_{k+p}P\wedge Q,
 \qquad
 P\vee Q=\mclass{0}{0_k}\oplus\mclass{Y}{1_p}.
\]
Indeed, the first inclusion follows from $T(x,0)=(Xx,0)$; the
second equality follows by adding arbitrary first coordinates
to the vectors $(Xx+Zy,Yy)$.

The valuation identity now gives
\[
 f(T)+kc=v(P)+v(Q)=v(P\wedge Q)+v(P\vee Q)
       \geq f(X)+kc+f(Y).
\]
Now we have $f(T)\geq f\left((X\oplus Y)\right)$. 
Finally, $\mclass{1_R}{1_R}=\mclass{0}{0_1}$, so $f(1_R)=c$.
\end{proof}

\begin{lemma}\label{lem:positive-block}
Let $v=(v_n)_{n\geq1}$ and $w=(w_n)_{n\geq1}$ be increasing
additive families of valuations vanishing at the origin. For each
$n$, let $u_n=v_n-w_n$ and let $u_n^+$ be its positive part on
$\mathcal L_n(R)$. Then
$u^+=(u_n^+)_{n\geq1}$ is an increasing additive family of valuations
vanishing at the origin.
\end{lemma}

\begin{proof}
For every chain $\mclass{0}{1_n}=P_0\leq_n\cdots\leq_n P_r=\mclass{0}{0_n}$,
\[
 \sum_{j=1}^r|u_n(P_j)-u_n(P_{j-1})|
 \leq v_n(\mclass{0}{0_n})+w_n(\mclass{0}{0_n}).
\]
Thus each $u_n$ has bounded variation. By Theorem~\ref{lem:birkhoff},
$u_n^+$ is increasing and vanishes at the origin. It remains
to prove additivity on direct sums.

Fix $P\in\mathcal L_n(R)$ and $Q\in\mathcal L_m(R)$. Define the maps
\[
\begin{aligned}
 \Phi\colon \{S\in\mathcal L_n(R)\mid S\leq_n P\}
 &\longrightarrow
 \{H\in\mathcal L_{n+m}(R)\mid H\leq_{n+m}P\oplus\mclass{0}{1_m}\},\\
 S&\longmapsto S\oplus\mclass{0}{1_m},\\[3pt]
 \Psi\colon \{T\in\mathcal L_m(R)\mid T\leq_m Q\}
 &\longrightarrow
 \{H\in\mathcal L_{n+m}(R)\mid H\leq_{n+m}\mclass{0}{1_n}\oplus Q\},\\
 T&\longmapsto \mclass{0}{1_n}\oplus T.
\end{aligned}
\]
Both maps are order-preserving and order-reflecting.
To determine their inverses, write
$H=\mclass{B}{\begin{pmatrix}C&D\end{pmatrix}}$, where $C$ has $n$
columns and $D$ has $m$ columns. For every left $R$-module $M$,
$(x,z)\in H(M)$ if and only if $By=Cx+Dz$ for some $y$.

If $H\leq_{n+m}P\oplus\mclass{0}{1_m}$, then
$H(M)\subseteq P(M)\times\{0\}$. Setting $z=0$, we obtain
\[
 H(M)=\mpair{B}{C}(M)\times\{0\},\qquad \mclass{B}{C}\leq_n P.
\]
Similarly, if $H\leq_{n+m}\mclass{0}{1_n}\oplus Q$, then
$H(M)\subseteq\{0\}\times Q(M)$. Setting $x=0$, we obtain
\[
 H(M)=\{0\}\times\mpair{B}{D}(M),\qquad \mclass{B}{D}\leq_m Q.
\]
Since these equalities hold for every $M$, the two maps are order
isomorphisms with inverses $H\mapsto\mclass{B}{C}$ and
$H\mapsto\mclass{B}{D}$ on their respective codomains.

Each map sends the least element to $\mclass{0}{1_{n+m}}$.
Since $u=v-w$ is additive and vanishes at the origin,
\[
 u_{n+m}(\Phi(S))=u_n(S),\qquad u_{n+m}(\Psi(T))=u_m(T).
\]
Thus the finite chains  correspond
under these maps, with the same increments on corresponding chains. Hence
\[
 u_{n+m}^+(P\oplus\mclass{0}{1_m})=u_n^+(P),\qquad
 u_{n+m}^+(\mclass{0}{1_n}\oplus Q)=u_m^+(Q).
\]
Since the two coordinate classes satisfy
\[
 \begin{aligned}
 (P\oplus\mclass{0}{1_m})\wedge(\mclass{0}{1_n}\oplus Q)
     &=\mclass{0}{1_{n+m}},\\
 (P\oplus\mclass{0}{1_m})\vee(\mclass{0}{1_n}\oplus Q)
     &=P\oplus Q.
 \end{aligned}
\]
The valuation identity for $u_{n+m}^+$ now gives
\begin{align*}
 u_{n+m}^+(P\oplus Q)
 & = u_{n+m}^+(P\oplus Q) +u_{n+m}^+( \mclass{0}{1_{n+m}} )   \\
 & = u_{n+m}^+(P\oplus\mclass{0}{1_m})+u_{n+m}^+(\mclass{0}{1_n}\oplus Q)  \\
   & =u_n^+(P)+u_m^+(Q).
\end{align*}
\end{proof}

\begin{proof}[Proof of Theorem~\ref{thm:ring}]
Let $E=\mathcal R(R)-\mathcal R(R)$. For $f\in E$, define $v_f$ as in Lemma~\ref{lem:rank-valuation}. This assignment is
linear. If $f=\rho-\sigma$ with $\rho,\sigma\in\mathcal R(R)$,
then $v_f=v_\rho-v_\sigma$ is a valuation family of bounded variation.
Moreover,
\[
 v_f(\mclass{X}{1_k})=f(X)\qquad(X\in M_{k,n}(R)).
\]
By Lemma~\ref{lem:positive-block}, $(v_f)^+$ is an increasing
additive family of valuations vanishing at the origin. Now apply  Corollary~\ref{lem:valuation-matrix}, define
\[
 h(X)=(v_f)^+(\mclass{X}{1_k})\qquad(X\in M_{k,n}(R)).
\]
Then $h\in\mathcal R(R)$.
Since $(v_f)^+-v_f$ is increasing, additive, and vanishes at the origin,
the same corollary gives
$h-f\in\mathcal R(R)$. Thus $h\geq0,f$.

Suppose
$$
g\geq0,f\quad {\rm in}\quad \mathcal R(R)-\mathcal R(R).
$$
 Then
$g,g-f\in\mathcal R(R)$, so $v_g$ and $v_g-v_f$ are increasing
by Lemma~\ref{lem:rank-valuation}. The least-upper-bound property
in Theorem~\ref{lem:birkhoff} gives $v_g\geq (v_f)^+$. Apply Corollary~\ref{lem:valuation-matrix} to $v_g-(v_f)^+$ to obtain
$g-h\in\mathcal R(R)$. Therefore $h=f\vee0$. Let $f^+=h$.
In particular,
\[
 f^+(X)=(v_f)^+(\mclass{X}{1_k}).
\]
Every element has a positive part.

Let $x,y\in E$. Since
\[
 (x-y)^+\geq x-y,\qquad (x-y)^+\geq0,
\]
the element $y+(x-y)^+$ is an upper bound of $x$ and $y$.
If $z\in E$ is any other upper bound of $x$ and $y$, then
\[
 z-y\geq x-y,\qquad z-y\geq0.
\]
Since $(x-y)^+=(x-y)\vee 0$, we obtain $(x-y)^+\leq z-y$, and hence
$y+(x-y)^+\leq z$. Thus
\[
 x\vee y=y+(x-y)^+.
\]
Similarly, the greatest lower bound also exists
and is given by
\[
 x\wedge y=-\bigl((-x)\vee(-y)\bigr).
\]
Therefore $E$ is a vector lattice.

If $\mathbb P(R)$ is nonempty, it is a compact base of
$\mathcal R(R)$. Since evaluation at $1_R$ is continuous
and strictly positive on the nonzero cone,
Proposition~\ref{prop:simplex} implies that $\mathbb P(R)$
is a Choquet simplex.
\end{proof}
\section{The simplex of dimension functions}\label{sec:refinement}
Throughout this section, $A$ is a unital $C^*$-algebra.
We prove Theorem~\ref{thm:dimension-simplex} by identifying the
dimension functions as a hereditary subcone of $\mathcal R(A)$.
We first recall the canonical extension to rectangular matrices.

For $X\in M_{p,q}(A)$, define the extension of a dimension function $d$ by
\[
 d(X)=d(|X|),\qquad |X|=(X^*X)^{1/2}.
\]
By Proposition~\ref{prop:cuntz-basic}, the same value is given by
$d(X^*X)$, $d(XX^*)$, or $d(|X^*|)$.

\begin{lemma}[{\cite[Proposition~5.2]{HL}}]\label{lem:matrix}
The canonical extension of a dimension function is a Sylvester
matrix rank function.
\end{lemma}

The cited result is stated for normalized functions. The present
form follows by scaling, since $d(1_A)=0$ implies $d=0$.
By \cite[Lemma~5.1]{HL},
\[
 X\preceq_{\mathrm M}Y\quad\Longrightarrow\quad X^*X\precsim Y^*Y.
\]
The converse fails in general \cite[Examples~5.4 and~5.5]{HL}.

The next lemma is the step that uses the factorization in
Proposition~\ref{prop:factor}. It shows that additive summands
of a dimension function preserve Cuntz comparison under the
stated product inequality.

\begin{lemma}\label{lem:summands}
Let $d$ be a dimension function, extended to rectangular
matrices as in Lemma~\ref{lem:matrix}. Let $f,g$ be functions from rectangular matrices over $A$ to
$[0,\infty)$ such that $f+g=d$. Suppose that, for $h=f,g$,
\[
 h(0)=0,\qquad
 h\begin{pmatrix}X&0\\0&Y\end{pmatrix}=h(X)+h(Y),\qquad
 h(XY)\leq\min\{h(X),h(Y)\}.
\]
Then the restrictions of $f$ and $g$ to positive matrices are
dimension functions, and $f,g$ are their canonical rectangular
extensions.
\end{lemma}

\begin{proof}
The product inequality gives
$f(X^*X)\leq f(X)$ and $g(X^*X)\leq g(X)$.
Since
\[
 f(X^*X)+g(X^*X)=d(X^*X)=d(X)=f(X)+g(X),
\]
both inequalities are equalities. The same argument applied to
$XX^*$ gives
\[
 f(X^*X)=f(X)=f(XX^*),\qquad
 g(X^*X)=g(X)=g(XX^*).
\]

Let $b\geq0$ and $c\in\Her(b)_+$. By
Proposition~\ref{prop:hereditary-basic}(1), $c\sim bcb$, so $d(bcb)=d(c)$.
On the other hand, the product inequality gives
$f(bcb)\leq f(c)$ and $g(bcb)\leq g(c)$.
Equality of their sums therefore implies
\[
 f(c)=f(bcb)\leq f(b),\qquad
 g(c)=g(bcb)\leq g(b).
\]

Now let $a\in M_p(A)_+$ and $b\in M_q(A)_+$ satisfy
$a\precsim b$. For every $n\geq1$, Proposition~\ref{prop:factor}
gives $T\in M_{(n+1)q,np}(A)$ such that
\[
 T^*T=\begin{pmatrix}a&&0\\&\ddots&\\0&&a\end{pmatrix},\qquad
 TT^*\in\Her\!\left(
 \begin{pmatrix}b&&0\\&\ddots&\\0&&b\end{pmatrix}\right).
\]
There are $n$ copies of $a$ and $n+1$ copies of $b$.
By the preceding observations,
\[
 nf(a)=f(T^*T)=f(TT^*)\leq(n+1)f(b).
\]
Dividing by $n$ and letting $n\to\infty$ gives $f(a)\leq f(b)$.
The same argument gives $g(a)\leq g(b)$. Thus both restrictions
are dimension functions, and the first part of the proof shows
that $f$ and $g$ are their canonical extensions.
\end{proof}

\begin{theorem}\label{thm:main}
The canonical extension
\[
 \iota:\DF(A)\longrightarrow\mathbb P(A),\qquad \iota(d)(X)=d(|X|),
\]
is an affine homeomorphism onto a closed face. The canonical extensions
of all dimension functions form the positive cone of a vector sublattice
of $\mathcal R(A)-\mathcal R(A)$.
\end{theorem}

\begin{proof}
By Lemma~\ref{lem:matrix}, $\iota$ takes values in $\mathbb P(A)$.
It is affine and continuous for pointwise convergence. Restriction
to positive matrices is its inverse on the image. Thus it is an
affine homeomorphism onto its image, which is closed because $\DF(A)$ is compact.

Suppose $\iota(d)=t\rho+(1-t)\sigma$, where $0<t<1$ and
$\rho,\sigma\in\mathbb P(A)$. Lemma~\ref{lem:summands}, applied
to $t\rho$ and $(1-t)\sigma$, shows that both restrictions are
dimension functions and that the summands are their canonical
extensions. Dividing by $t$ and $1-t$, respectively, gives
$\rho,\sigma\in \iota(\DF(A))$. Hence the image is a face.

Identify dimension functions with their canonical extensions, and
let $C\subseteq\mathcal R(A)$ be their cone. Lemma~\ref{lem:summands}
shows that $C$ is hereditary in the sense of \cite[\S1.1]{AG}:
\[
 0\leq f\leq d,\quad d\in C\quad\Longrightarrow\quad f\in C.
\]
Here the order is the algebraic order on $\mathcal R(A)$.
In particular, if $d,e\in C$ and $e-d\in\mathcal R(A)$, then
$0\leq e-d\leq e$ implies $e-d\in C$. Thus the algebraic order
on $C$ agrees with the restricted order, and
$(C-C)\cap\mathcal R(A)=C$.

By Theorem~\ref{thm:ring}, $\mathcal R(A)-\mathcal R(A)$ is a
vector lattice. For $d,e\in C$ we have
$0\leq(d-e)^+\leq d$, so $(d-e)^+\in C$. Hence $C-C$ is a
vector sublattice with positive cone $C$.
\end{proof}

\begin{proof}[Proof of Theorem~\ref{thm:dimension-simplex}]
Suppose that $\DF(A)$ is nonempty. It is the compact base of
the cone of dimension functions defined by $d(1_A)=1$.
Evaluation at $1_A$ is continuous and strictly positive on the
nonzero cone. Theorem~\ref{thm:main} and
Proposition~\ref{prop:simplex} therefore show that $\DF(A)$ is a
Choquet simplex.
\end{proof}

Finally, we provide a characterization of when $\mathrm{DF}(A)$ is empty.
\begin{proposition}\label{prop:empty-df}
Let $A$ be a nonzero unital $C^*$-algebra. Then $\DF(A)$ is empty
if and only if $1_{M_n(A)}$ is properly infinite for some $n\geq1$.
\end{proposition}

\begin{proof}
Let $u=[1_A]$. For every $a\in M_m(A)_+$, $[a]\leq mu.$  Follow the notation of \cite{OPR}, let $S(W(A),u)$ be the set of
additive order-preserving maps $f:W(A)\to[0,\infty]$ such that
$f(u)=1$. Every such map satisfies
\[
 0\leq f([a])\leq m \qquad(a\in M_m(A)_+),
\]
and therefore takes finite values. The assignments
\[
 d_f(a)=f([a]),\qquad f_d([a])=d(a)
\]
identify $S(W(A),u)$ with $\DF(A)$.

In \cite[Proposition~2.1]{OPR}, take $x=y=u$.
The condition $x\leq ky$ for some $k\geq1$ holds with $k=1$.
The remaining condition is
\[
 f(u)<f(u)\qquad\text{for every }f\in S(W(A),u).
\]
Since $f(u)=1$, this condition holds precisely when $S(W(A),u)$
is empty. Hence
\[
 \DF(A)=\varnothing
 \quad\Longleftrightarrow\quad
 (n+1)u\leq nu \quad\text{for some }n\geq1.
\]
Furthermore, $(n+1)u \le nu$ holds if and only if $2nu \le nu$. For the projection $p = 1_{M_n(A)}$, we obtain $[p] = nu$, so that $p \oplus p \precsim p$, which means that $p$ is a properly infinite projection.
\end{proof}

Theorem~\ref{thm:dimension-simplex} now shows that, for nonzero unital
$A$, $\DF(A)$ is a nonempty Choquet simplex precisely when no $M_n(A)$
has a properly infinite unit.

\section*{AI statement}
The authors used GPT 6 Astra during the development of this work for searching necessary references, examining mathematical arguments, and language editing.  The authors have checked and substantially revised all mathematical arguments and takes full responsibility for the final version.

\end{document}